\documentclass[11pt]{article}
\usepackage{graphicx}
\newtheorem{theorem}{Theorem}
\newtheorem{lemma}[theorem]{Lemma}

\newtheorem{proposition}[theorem]{Proposition}

\newtheorem{definition}[theorem]{Definition}

\newtheorem{notation}[theorem]{Notation}
\newtheorem{remark}[theorem]{Remark}

\begin{document}

\title{Sets of double and  triple weights of trees}
\author{  Elena Rubei }
\date{\hspace*{1cm}}
\maketitle

\vspace{-1.1cm}

{\small 
{\bf Address:}
Dipartimento di Matematica ``U. Dini'', 
viale Morgagni 67/A,
50134  Firenze, Italia 

{\bf
E-mail address:}
 rubei@math.unifi.it}

\bigskip

\smallskip

\def\thefootnote{}
\footnotetext{ \hspace*{-0.36cm}

}

{\small {\bf Abstract.} 
Let $T$ be  a weighted tree   with $n$ leaves numbered by the set 
$\{1,...,n\}$.
Let $D_{i,j}(T)$ be the distance between the leaves $i$ and $j$. 
Let $D_{i,j,k}(T)= \frac{1}{2}\left( D_{i,j}(T) + D_{j,k}(T) +D_{i,k}(T)
 \right)$.
We will call such numbers ``triple weights'' of the tree.
 
In this paper, we give a characterization, different from
 the previous ones, for sets indexed by 
$2$-subsets of a $n$-set to be double weights of a tree.

By using the same ideas,
we find also necessary and sufficient conditions for a set of 
 real numbers  indexed by  $3$-subsets of an $n$-set 
to be the set of the triple weights of a tree with $n$ leaves.

Besides we propose a slight modification of Saitou-Nei's Neighbour-Joining 
algorithm to reconstruct trees from the data $D_{i,j} $. }


\bigskip

\section{Introduction}

Consider a positive-weighted tree $T$ (that is a tree such that every 
edge is endowed with a positive real number, which we  call the 
length of the edge) with $n$ leaves numbered by the set 
$\{1,...,n\}$. 
Let $D_{i,j}(T)$ be the sum of the lengths of the edges of the shortest 
path connecting $i$ and $j$. We call such number the ``distance''
 between the leaves $i$ and $j$ or the ``double weight'' for $i$ and $j$. 

In 1971 Buneman characterized the metrics on finite sets coming from a tree:

\begin{theorem} {\bf (Buneman)} A metric $(D_{i,j})$ 
on $\{1,...,n\}$ is the metric 
induced
 by a positive-weighted tree if and only if for all $i,j,k,h  \in \{1,...,n\}$
the maximum of $\{D_{i,j} + D_{k,h},D_{i,k} + D_{j,h},D_{i,h} + D_{k,j}
 \}$ 
is attained at least twice. 
\end{theorem}

The problem of reconstructing trees from data involving the distances  
between the leaves has several applications, such as internet tomography and 
phylogenetics;
 evolution of species can be represented by trees and, given
 distances between genetic  sequences of some species, one 
can  try to reconstruct the  evolution tree from these distances.
Some algorithms to reconstruct trees from the data $\{D_{i,j}\}$ have been 
proposed. Among them is neighbour-joining method, invented by
 Saitou and Nei in 1987 (see \cite{NS},
\cite{SK} and \cite{PSt2}).

For any weighted  tree $T$, 
let now  $$D_{i,j,k}(T)= \frac{1}{2}( D_{i,j}(T) + D_{j,k}(T) +D_{i,k}(T)),$$
that is the sum of the lengths of the edges of the minimal subtree 
with $i,j,k$ as set of leaves. 
We  call such numbers ``triple weights'' of the tree.
More generally define the $k$-weights of the tree
 $D_{i_1,....,i_k}(T)$ as the sum of the lengths of the edges of the minimal 
subtree connecting $i_1$,....,$i_k$.

In 2004, Pachter and Speyer proved the following theorem.

\begin{theorem} {\bf (Pachter-Speyer)}. Let $ k ,n  \in {\bf N}$ with
$ n \geq 2k-1$ and $ k \geq 3$.  
A positive-weighted tree
 $T$ with $n$ leaves $1,...,n$ and no vertices of degree 
2
is determined by the values $D_I$ where $ I $ varies in the $k$-subsets 
of $\{1,...,n\}$.
\end{theorem}

It can be interesting to characterize the sets of real numbers 
which are sets of $k$-weights of a tree, for instance triple weights of a
tree; in fact 
(I quote Speyer and Sturmfels's  paper \cite{SS2})
``it can be more reliable statistically
to estimate the triple weights $D_{i,j,k}$ rather than the pairwise
distances $D_{i,j}$''.  We refer to \cite{SS2} and above all to 
\cite{PS} for an analysis of this and the references.

In \cite{BC} Bocci and Cools 
give a description of $k$-dissimiliarity maps of a trees
and generalize Buneman's result for
sets of real numbers $\{D_{i,j,l}\}$ indexed by $3$-subsets of $\{1,...,n\}$ 
coming from sets $\{D_{i,j}\}$.

In this paper we give a characterization (different from 
Buneman's one) for sets indexed by 
$2$-subsets of a $n$-set to be double weights of a tree with $n$ leaves
(see Theorem \ref{2weights}) (please note that in this paper a weight 
is not necessarily positive).

By using the same ideas, 
we find also necessary and sufficient conditions for 
a set of real numbers  indexed by  $3$-subsets of an $n$-set 
to be the set of the triple weights of a tree  (see 
Theorem \ref{3weightsbis}).

Finally, by using the characterization of neighbours we used to deduce the 
above  theorems,   we propose a slight modification of neighbour-joining
algorithm.

\section{The main theorems}


\begin{definition} \label{bell}
A {\bf 2-cherry}  $B$ in a tree $T$ is 
a subtree $B$ with two leaves such that only one of the inner vertices is not 
bivalent;
we call this vertex  ``{\bf stalk}'' of the bell
and we say that the two leaves  are {\bf neighbours}. 
We call the path from a leaf of a bell to its stalk 
``{\bf twig}'' of this leaf.

A {\bf cherry} is a union of 2-cherries with the same stalk.

\end{definition}

\begin{center}
\includegraphics[scale=0.45]{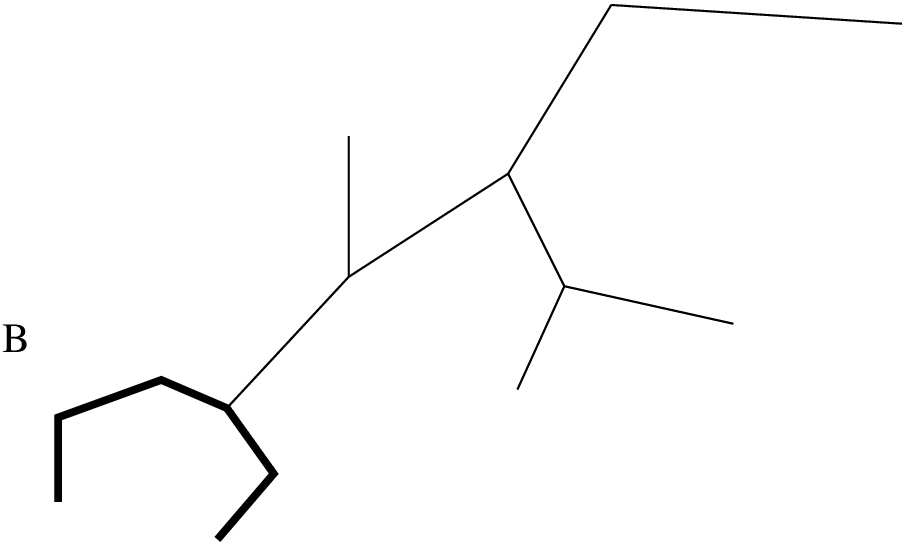}
\end{center}

\begin{notation}
$\bullet$ For every $n \in {\bf N}$ let $ [n] = \{1,....,n\}$; besides, 
if $ k \leq n$, denote  the set of the $k$-subsets of $[n]$ by $[n]_k$. 

\smallskip

$\bullet$ 
 For any  set $\{ D_{\{i_1,...., i_k\}}\}$ of real numbers indexed by 
the $k$-subsets  $ \{i_1,...., i_k\} $  of $ [n]$,
we  denote $D_{\{i_1,..., i_k\}}$ by $ D_{i_1,..., i_k}$ 
for any order of $i_1,...., i_k$.

\smallskip

$\bullet$ A {\bf weighted tree} is a tree such that every edge is 
endowed with a real number called weight or length of the edge. If 
the weights are positive we say that the tree is {\bf positive-weighted}.
Please note that in other papers ``weighted'' means positive-weighted.

\end{notation}

\begin{remark} \label{dt}
Let 
$\{D_{i,j}\}_{\{i,j\} \subset [n]_2}$
be a set of real numbers 
and let  $D_{i,j,k}:= \frac{1}{2}( D_{i,j} + D_{j,k} +D_{i,k})$
Then
$$ D_{i,j}= \frac{2}{3}
(D_{i,j,r}+ D_{i,j,s}+ D_{i,j,u}+ D_{r,s,u}) -
\frac{1}{3}(
D_{i,r,s}+D_{i,r,u}+D_{i,s,u} +D_{j,r,s}+D_{j,r,u}+D_{j,s,u})$$ 
for any $r,s,u$. 
\end{remark}

\begin{notation}
For any  set $\{ D_{i_1,...., i_k}\}_{\{i_1,...., i_k\} \subset  [n]_k}$
 of real numbers   and for any 
$\alpha, \alpha' \in [n]$,
we  define $\ast^{\alpha, \alpha'} $ the following condition:
\begin{center}
$  D_{\alpha, \gamma_1,...., \gamma_{k-1}}  - D_{\alpha',
 \gamma_1,...., \gamma_{k-1} } $
doesn't depend on $\gamma_1,...., \gamma_{k-1}
\in [n] - \{\alpha , \alpha'\}$
\end{center}
\end{notation}

\begin{proposition} \label{starbell}
Let $T$ be a positive-weighted tree with leaves  $1,....,n$ with 
$n \geq 2k-1$.  Let $ \alpha , \alpha' \in [n]$.
Then $  \ast^{\alpha, \alpha'} $ holds for $ \{ D_{i_1,...,i_k} (T) \}_{\{
i_1,...
,i_k\} \in [n]_k}$ if and only if $ \alpha$  and $ \alpha' $ form a cherry.
\end{proposition}

{\it Proof.}

$\Leftarrow$ Obvious.

$\Rightarrow$ We say that an inner 
 vertex is a node if it is not bivalent.
Let $ \alpha , \alpha'$ be such that  $  \ast^{\alpha, \alpha'} $ holds.
Let us suppose that they  are not neighbours. 
Then on the path from $\alpha $ to $\alpha' $  there are at least two nodes.
Let $x, y$ be two nodes on the path from $\alpha$ to $ \alpha' $ such that 
there is no node in the path from $x$ to $y$ (suppose $x$ nearer to $\alpha$ 
and $y$ nearer to $\alpha'$). 
For every $ \delta \in [n]$, let $ \overline{\delta}$ be the 
node on the path from $\alpha $ to $ \alpha'$ such that 
$$ path(\alpha, \alpha') \cap path(\alpha, \delta) = 
path(\alpha , \overline{\delta})$$

We can divide $[n]$ into two disjoint subsets:

$X =\{\delta \in [n] \; |\; \overline{\delta}\; is\; between \;\alpha \; 
and \;x\}$

$Y =\{\delta  \in [n] \; |\; \overline{\delta}\; is\; between \;y \; 
and \;\alpha'\}$

Since $ n \geq 2k-1$, one of the two subsets must contain at least $k-1$ 
elements $\gamma_1,..., \gamma_{k-1}$ besides one among  $ \alpha$ or 
$\alpha'$; we can suppose it is $X$; the other one, that is $Y$,  must 
contain another element  $\eta $ besides $\alpha'$ (since $y$ 
is a node).
Up to changing the names of $\gamma_1,..,\gamma_{k-1}$  (and 
correspondingly $\overline{\gamma_1},...,\overline{\gamma_{k-1}})$,
we can suppose 
$$length\;path (\alpha, \gamma_1) \leq length\;path (\alpha, \gamma_2) 
\leq..... \leq length\;path (\alpha, \gamma_{k-1} )$$

\begin{center}
\includegraphics[scale=0.45]{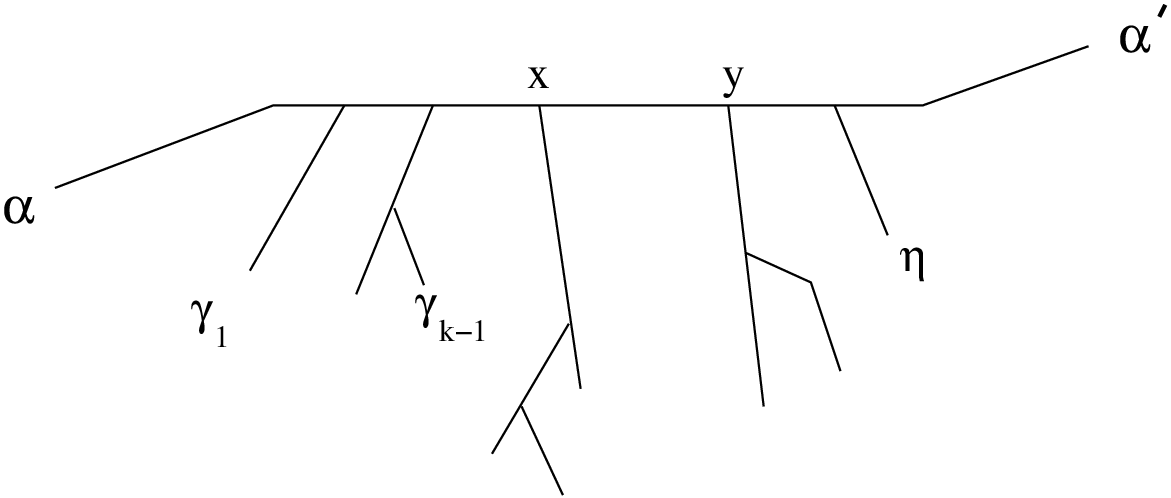}
\end{center}

Then we have 
$$  D_{\alpha, \gamma_1,...., \gamma_{k-1}}  - D_{\alpha',
 \gamma_1,...., \gamma_{k-1} } 
= length\; path (\alpha,\overline{\gamma_{k-1}} ) -
length\; path (\alpha',\overline{\gamma_{1}} ) <$$
$$ <length\; path (\alpha,\overline{\eta}) -
length\; path (\alpha',\overline{\gamma_{2}} )=
 D_{\alpha,  \eta,\gamma_2,...., \gamma_{k-1}}  - D_{\alpha',
 \eta, \gamma_2,...., \gamma_{k-1} } $$
contradicting the assumption $\ast^{\alpha, \alpha'}$.
\hfill \framebox(7,7)

\begin{remark}
The trees with $4$ or $5$ leaves are homeomorphic  to the 
following trees (with the nodes possibly collapsed):
\end{remark}

\begin{center}
\includegraphics[scale=0.45]{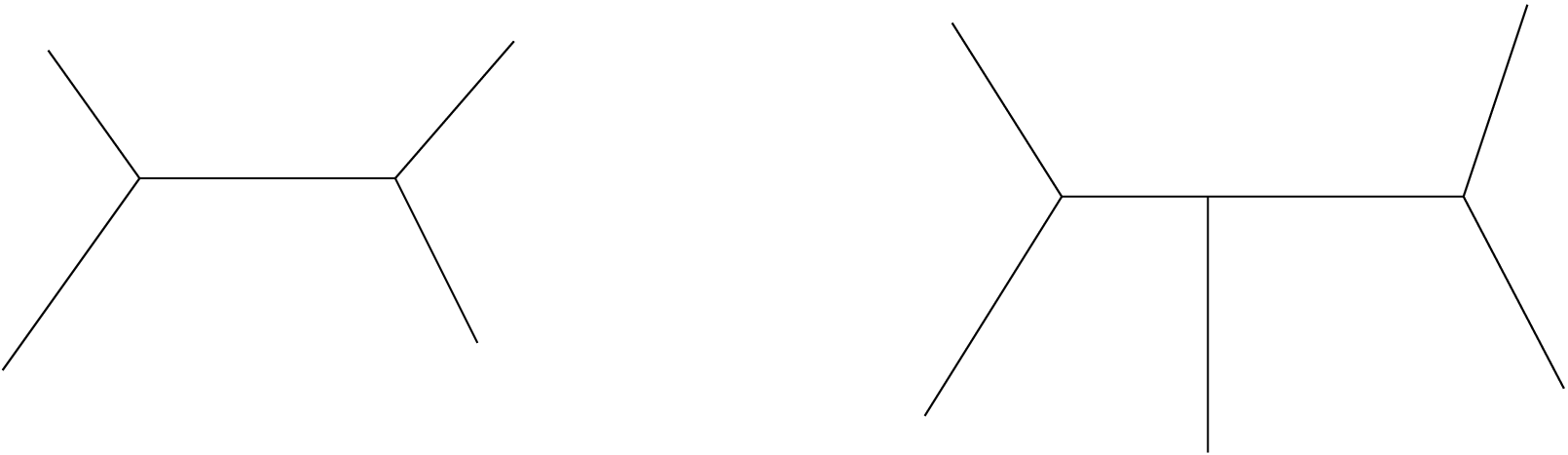}
\end{center}

\begin{proposition} \label{5}
1) Let $\{ D_{i,j}\}_{\{i,j\} \in [4]_2}$ be a set of real
 numbers.
There exists a weighted tree $T$ such that  $ D_{i,j} =   D_{i,j} (T)$ 
if and only if there exist 
$\alpha, \alpha', $ distinct in $ [4]$ 
such that 
$\ast^{\alpha, \alpha'}$ holds (observe that then also
$\ast^{\beta, \beta'}  $ holds where $ \{\beta, \beta'\} = [4] - \{\alpha,
\alpha'\}$).

2) Let $\{ D_{i,j,k}\}_{\{i,j,k\} \in [5]_3}$ be a set of real
 numbers.
There exists a weighted tree $T$ such that  $ D_{i,j,k} =   D_{i,j,k} (T)$ 
if and only if there exist 
$\alpha, \alpha', \beta, \beta',
 $ distinct in $ [5]$ 
such that 
$\ast^{\alpha, \alpha'}$  and 
$\ast^{\beta, \beta'}  $ hold.
\end{proposition}

{\em Proof.} We prove 1 (2 is analogous).
Let $T$ be a weighted tree with $4$ leaves. By the previous remark 
it is homeomorphic to the following tree and we call the leaves and the
 weights  as
  in the figure:

\begin{center}
\includegraphics[scale=0.45]{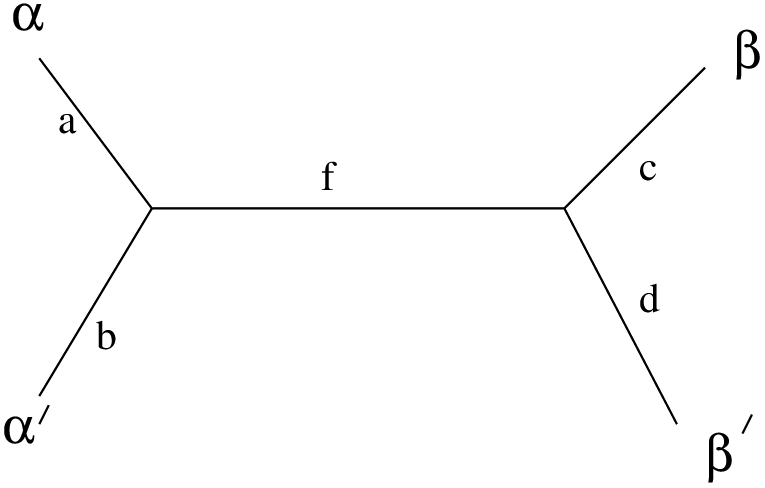}
\end{center}

Obviously we have:

$ D_{\alpha , \alpha'} = a+b$ 

$ D_{\alpha , \beta } = f +a+c$ 

$ D_{\alpha , \beta'} = f +a+d$ 

$ D_{ \beta, \beta' } = c+d$ 

$ D_{\alpha' , \beta } = f +b+c$ 

$ D_{\alpha' , \beta'} = f +b+d$

Let $$ S := 
\left\{ {\scriptsize \left(  
\begin{array}{c} 
 a+b \\
 f +a+c \\
  f +a+d \\
 c+d  \\
  f +b+c  \\
 f +b+d 
\end{array} \right)}
| \; a,b,c,d,f \in {\bf R}\right\}$$
One can easily see that 
$$ S =\left\{ 
{\scriptsize \left(  
\begin{array}{c} 
 D_{\alpha , \alpha'} \\
 D_{\alpha , \beta } \\
 D_{\alpha , \beta'} \\
 D_{ \beta, \beta' } \\
 D_{\alpha' , \beta } \\
 D_{\alpha' , \beta'} 
\end{array} \right)} 
\in {\bf R}^{6} | \; \; \ast^{\alpha, \alpha'} \; (and \;
\ast^{\beta, \beta'}) \; hold(s) 
 \right\}$$
\hfill \framebox(7,7)

\begin{notation} \label{pseudobell}
For every set  $\{ D_{i_1,...,i_k}\}_{\{i_1,...., i_k\} \in [n]_k} $ 
of  real numbers 
we say that $\underline{\alpha} = \{\alpha_1,....,\alpha_r\} \subset  
 [n]$ is a {\bf ($r$-)pseudocherry}  
if $ \ast^{\alpha_i, \alpha_j}$ holds
for all $ i,j $. 

We say that $\underline{\alpha}$ is a   {\bf complete pseudocherry}
if $\not \exists \beta \in [n]- \underline{\alpha}$
 such that $ \ast^{\beta, \alpha_i}$  holds for all 
$ i $. 
\end{notation}

Now we are ready to state the theorems characterizing 
the sets of real numbers indexed by $2$-subsets (or $3$-subsets) of $[n]$ 
which come from a tree.
Shortly speaking such a set comes from a tree if and only if in 
$[n]$ there are at least two pseudocherries and if we substitute
every pseudocherry with a point, the same condition holds for the new set
 and so on.

\begin{theorem}  \label{2weights} Let $n \geq 4$. 
Let $\{ D_{i,j}\}_{\{i,j \} \in [n]_2}$ be a set of  real numbers. 
There exists a tree $R$  such that  $D_{i,j} = D_{i,j}(R) $ 
for all $ i,j$  if and only if  the following conditions hold:

- There exists at least one 2-pseudocherry in $ [n]
\;^{1}$. 
\def\thefootnote{1}
\footnotetext{ We could state the theorem also by saying ``There exist
at least two disjoint 2-pseudocherries in $[n]$''.}

Define $L$ as the set obtained from  $[n]$ by substituting, for any 
complete 
pseudocherry $\underline{\alpha}$, all the elements of the pseudocherry with
a new   leaf  $  z_{\underline{\alpha}}  $ and define for all $\alpha \in 
\underline{\alpha}$
$$a_{\alpha}:=\frac{1}{2} ( D_{ \alpha , 
\alpha'}+ D_{ \alpha , x} -D_{ \alpha',x})$$
for any  $\alpha' \in \underline{\alpha}$, $x\in [n] -\{\alpha, \alpha'\}$.

For any  $\underline{\alpha}$ pseudocherry, and possibly  
$\underline{\beta}$, pseudocherry disjoint from
$\underline{\alpha}$,  define
 $$ D_{ z_{\underline{\alpha}},j} :=D_{ \alpha,j} - a_{\alpha} \;\; \forall 
\alpha \in \underline{\alpha}$$
 $$ D_{ z_{\underline{\alpha}},z_{\underline{\beta}} } :=
D_{ \alpha, \beta} - a_{\alpha} - a_{\beta}  \;\; \forall 
\alpha \in \underline{\alpha} , \;\;
\beta \in \underline{\beta}$$

-Then the same condition holds  for $D_{i,j}$ with  $i,j \in L$ 
instead of  $ [n]$ and so on up to  a 4-set (if at the last step 
you get  a set with less than 4 elements, take off less pseudocherries).

\medskip

Obviously $R$ i is positive-weighted if and only if
$a_{\alpha}$ is positive for any $\alpha \in [n]$ and also for any
$\alpha \in L$ and so on up to a set of $3$-elements (that is $R$ is 
positive-weighted 
if and only if $a_{\alpha} >0$  for all $\alpha $ in the sets 
we get by ``deleting'' pseudocherries until we get a set of $3$ elements).

\end{theorem}

{\it Proof.}

$\Rightarrow$ Obvious.

\medskip

$\Leftarrow$ By induction on  $n$. 
The case  $n=4$ follows from Proposition \ref{5}.

Let us prove the induction step.

One can easily prove that 
the definition of $a_{\alpha}$ for $\alpha$
in a pseudocherry $\underline{\alpha}$ is a good definition, that is it 
doesn't depend either on $\alpha'$ or on  $x,y$ (because 
$\ast^{\alpha, \alpha'}$ holds). 
We have to prove  also that, for any  pseudocherry $\underline{\alpha}$,
the definition of $ D_{ z_{\underline{\alpha}},j} $,  
 is a good definition, 
in fact  we have defined it as  
$ D_{\alpha,j} -a_{\alpha} $
for any $x $ and  for some $\alpha \in \underline{\alpha}$, 
but we could define  it as  
$ D_{\alpha',j} -a_{ \alpha'}$
 for another  element 
$\alpha' \in \underline{\alpha}$. 
Obviously they are the same (choose $x=j$  in the expressions
of $a_{ \alpha}$ and $a_{ \alpha'}$).
Analogously we can prove that the definitions of 
$ D_{ z_{\underline{\alpha}},z_{\underline{\beta}}}$ 
is a   good definition.

By induction assumption we can suppose that there exists a tree 
$T$ such that $$  D_{i,j} =  D_{i,j}(T)$$ for any $i,j \in  L$. 
We define $R$ as the tree obtained from $T$ by adding, for every pseudocherry 
$ \underline{\alpha}$ of  $[n]$, 
 a cherry  with leaves the elements of the pseudocherry
and  stalk $z_{\underline{\alpha}}$ and $a_{ \alpha}$ as length of the twig
of $\alpha $  for any  $\alpha \in  \underline{\alpha}$.

\medskip

We have to show that  $ D_{i,j} (R) = D_{i,j} $ for any $i,j
 \in [n]$. If $ i,j $ are  elements of $L$ too, we already know that.
So we have to show that  $ D_{i,j} (R) = D_{i,j} $ for any $i,j
 \in [n]$ with at least one of $i,j$ not in $L$, that is with at least one
 of $i,j$ in a pseudocherry of $[n]$.

We can divide the proof into  three  cases:

\medskip

1) 
$D_{ \alpha , j} (R) =D_{ \alpha , j}  $ for any $j \in [n] \cap L$
 and   $\alpha$ in a pseudocherry $\underline{\alpha}$ in  $[n]$

2)
$D_{ \alpha , \beta } (R) =D_{ \alpha , \beta }  $ for any 
  $\alpha \in \underline{\alpha}$,
$\beta \in \underline{\beta}$ with  
 $ \underline{\alpha}, \underline{\beta}$ disjoint pseudocherries
in  $[n]$

3)
$D_{ \alpha , \alpha' } (R) = D_{\alpha , \alpha'} $ 
for any    $\alpha, \alpha'$ in a pseudocherry 
$\underline{\alpha}$ in $[n]$

\bigskip

1) $D_{ \alpha , j} (R) \stackrel{\footnotesize de\!f. \; o\!f\; R}{=}\; 
a_{\alpha} + D_{z_{ \underline{\alpha}} , j}(T) 
\stackrel{\footnotesize ind. \;ass. }{=}\;
a_{\alpha} + D_{z_{ \underline{\alpha}} , j}
 \stackrel{\footnotesize de\!f. \; o\!f\  D_{z_{ \underline{\alpha}} , j} 
}{=}\;  
a_{\alpha} + D_{ \alpha , j} -a_{\alpha} =D_{ \alpha , j} $

\bigskip

2)  $D_{ \alpha , \beta } (R) \stackrel{\footnotesize de\!f. \; o\!f\; R}{=}\;
 D_{z_{\underline{\alpha}}, z_{\underline{\beta}}}(T) + a_{\alpha} 
+ a_{\beta} 
\stackrel{\footnotesize ind. \;ass. }{=}\;
 D_{z_{\underline{\alpha}}, z_{\underline{\beta}}} + a_{\alpha} + a_{\beta} 
 \stackrel{\footnotesize de\!f. \; o\!f\  
D_{z_{ \underline{\alpha}},z_{\underline{\beta}}} }{=}\;  
D_{\alpha, \beta}- a_{\alpha} - a_{\beta} 
+ a_{\alpha} + a_{\beta} =
D_{ \alpha , \beta }  $

\bigskip

3) 
$D_{ \alpha , \alpha' } (R)  
\stackrel{\tiny de\!f. \; o\!f\; R}{=}\; 
a_{\alpha} + a_{\alpha'} 
\stackrel{\forall \;x \in [n]-\{\alpha, \alpha'\} }{=}$

\bigskip

$=\frac{1}{2} 
( D_{\alpha, \alpha'}  + D_{\alpha , x} -  D_{\alpha' , x} ) +
 \frac{1}{2} 
( D_{\alpha', \alpha}  + D_{\alpha', x} -  D_{\alpha , x} ) 
= D_{ \alpha , \alpha' }$

\hfill \framebox(7,7)

\bigskip

In \cite{BC} Bocci and Cools 
generalize Buneman's result for
sets of real numbers $\{D_{i,j,l}\}_{\{i,j,l \} \in [n]_3}$
``coming from'' sets $\{D_{i,j}\}_{\{i,j \} \in [n]_2}$. 
Here, by using the same ideas used for Theorem \ref{2weights}, we 
give another  characterization for these sets to be sets of weigths of a 
tree.

\begin{theorem} \label{3weightsbis} Let $ n \geq 5$.
Let $\{ D_{i,j, k}\}_{\{i,j,k\} \in [n]_3}$ be a set of  real numbers
coming from a set $\{D_{i,j}\}$  
(by the formula $D_{i,j,k} := \frac{1}{2} (D_{i,j} + D_{i,k} + D_{k,j})$). 
There exists a tree $R$  such that  $D_{i,j,k} = D_{i,j,k}(R) $ 
for all  $i,j,k$  if and only if  the following conditions hold:

- There exist at least two disjoint 2-pseudocherries in $ [n]
\;^{1}$. 
\def\thefootnote{1}
\footnotetext{ This means that there exist at least distinct 
$ \alpha, \alpha' \beta ,
 \beta' $ such that $\ast^{\alpha , \alpha'}$ and 
$ \ast^{\beta, \beta'}$ hold;
we are not excluding that also $  \ast^{\alpha, \beta}$ 
may hold, that is, that $\alpha, \alpha' \beta , \beta' $ 
may be in the same pseudobell} 

Define $L$ as the set obtained from  $[n]$ by substituting, for any 
complete 
pseudocherry $\underline{\alpha}$, all the elements of the pseudocherry with
a new   leaf  $  z_{\underline{\alpha}}  $ and define for all $ \alpha \in 
\underline{\alpha}$
$$a_{\alpha}:=\frac{1}{2} ( D_{ \alpha , 
\alpha'}+ D_{ \alpha , x,y} -D_{ \alpha',x,y })$$
for any  $\alpha' \in \underline{\alpha}$, $x,y \in [n] -\{\alpha, \alpha'\}$
(where 
$D_{ \alpha , \alpha'}$ is given in function of the $ D_{i,j,k}$ by Remark 
\ref{dt}).

For any  $\underline{\alpha}$ pseudocherry, and possibly  
$\underline{\beta}$ and
 $\underline{\gamma}$ pseudocherries, all disjoint, define
 $$ D_{ z_{\underline{\alpha}},j,k} :=D_{ \alpha,j,k} - a_{\alpha} 
\;\; \forall 
\alpha \in \underline{\alpha}$$
 $$ D_{ z_{\underline{\alpha}},z_{\underline{\beta}} ,j} :=
D_{ \alpha, \beta,j} - a_{\alpha} - a_{\beta}  \;\; \forall 
\alpha \in \underline{\alpha} , \;\;
\beta \in \underline{\beta}$$
$$ D_{ z_{\underline{\alpha}},z_{\underline{\beta}} , z_{\gamma}} :=D_{ \alpha, \beta, \gamma} -
 a_{\alpha} -
 a_{\beta} -
a_{\gamma} \;\; \forall 
\alpha \in \underline{\alpha},   
\;\;\beta \in \underline{\beta}, 
\;\;\gamma \in \underline{\gamma}$$

-Then the same condition holds 
for $D_{i,j,k}$ with  $i,j,k \in L$ 
instead of  $ [n]$ and so on up to  a 5-set (if at the last step 
you get  a set with less than 5 elements, take off less pseudocherries).

\medskip

Obviously $R$ 
 is positive-weighted if and only 
 if $a_{\alpha} >0$  for all $\alpha $ in the sets 
we get by ``deleting'' pseudocherries until we get a set $S$ of $5$ elements
and also for every $\alpha $ in all the sets of $4$ elements 
we get by deleting two neighbours of $S$. 

\end{theorem}

{\it Sketch of the proof.}
The proof is analogus to the previous one.

By induction assumption we can suppose that there exists a tree 
$T$ such that $ D_{i,j,k} =  D_{i,j,k}(T)$ for any $i,j,k \in  L$. 
We define $R$ as the tree obtained from $T$ by adding, for every pseudocherry 
$ \underline{\alpha}$ of  $[n]$, 
 a cherry  with leaves the elements of the pseudocherry
and  stalk $z_{\underline{\alpha}}$ and $a_{ \alpha}$ as length of the twig
of $\alpha $  for any  $\alpha \in  \underline{\alpha}$.  
We  have to show that  $ D_{i,j,k} (R) = D_{i,j,k} $ for any $i,j,k
 \in [n]$ with some of $i,j,k$ not in $L$, that is with some of $i,j,k$
in a pseudocherry of $[n]$.

We can distinguish two main cases:

A) the elements among $ i,j,k$ that are in some pseudocherries of $[n]$, are in 
disjoint pseudocherries of $[n]$

B) there are at least two among $i,j,k$ in the same pseudocherry of $[n]$.

\medskip

The proof of case A can be divided into three subcases:

A1) 
$D_{ \alpha , j,k} (R) =D_{ \alpha , j,k}  $ for any $j,k \in [n] \cap L$
 and   $\alpha$ in a pseudocherry $\underline{\alpha}$ in  $[n]$

A2)
$D_{ \alpha , \beta , j} (R) =D_{ \alpha , \beta ,j}  $ for any 
$j \in [n] \cap L$,  $\alpha \in \underline{\alpha}$,
$\beta \in \underline{\beta}$ with  
 $ \underline{\alpha}, \underline{\beta}$ disjoint pseudocherries
in  $[n]$

A3)
$D_{ \alpha , \beta , \gamma} (R) =D_{ \alpha , \beta ,\gamma} $ 
 for any  $\alpha \in \underline{\alpha}$,  
 $\beta \in \underline{\beta}$,
 $\gamma \in \underline{\gamma}$ 
with  
 $ \underline{\alpha}, \underline{\beta}$,  $ \underline{\gamma}$ 
 disjoint pseudocherries in $[n]$.

The proof of 
case B can be divided into:

B1)
$D_{ \alpha , \alpha' , j} (R) = D_{\alpha , \alpha', j} $ 
for any $j \in [n] \cap L$,   $\alpha, \alpha'$ in a pseudocherry 
$\underline{\alpha}$ in $[n]$

B2) 
$D_{ \alpha , \alpha' , \beta} (R) = 
D_{\alpha , \alpha', \beta} $ 
for any  $\alpha, \alpha' \in\underline{\alpha}$,
 $\beta \in \underline{\beta}$ 
with  
 $ \underline{\alpha}, \underline{\beta}$ disjoint pseudocherries
in $[n]$

B3)
$D_{ \alpha , \alpha' , \alpha''} (R) = 
D_{\alpha , \alpha', \alpha''} $ 
 for any  $\alpha, \alpha', \alpha''$ in a pseudocherry 
$\underline{\alpha}$ in $[n]$.

We omit the proofs of cases A and we prove only 
the (more complicated) cases  B.





\bigskip

B1) $D_{ \alpha , \alpha' , j} (R)  
\stackrel{\tiny de\!f. \; o\!f\; R}{=}\;
a_{\alpha}+ a_{\alpha'} + D_{z_{\underline{\alpha}} , j}(T) =$

\bigskip

$\stackrel{\tiny \forall\;x,y,u \in L, \;\; Rem. \ref{dt}}{=}\; 
 a_{\alpha}+ a_{{\alpha'}} + \frac{2}{3}[ D_{z_{\underline{\alpha}}, j,x}(T) 
+ D_{z_{\underline{\alpha}}, j,y}(T)  + D_{z_{\underline{\alpha}}, j,u}(T) + D_{x,y,u} (T)]
-\frac{1}{3} [ D_{z_{\underline{\alpha}}, x,y}(T) 
+ D_{z_{\underline{\alpha}}, y,u}(T)  + D_{z_{\underline{\alpha}}, x,u}(T) + D_{j,x,y} (T)
+ D_{j,y,u} (T) + D_{j,x,u} (T)]=$

\bigskip

$ \stackrel{\tiny ind. \;ass. }{=}\; 
 a_{\alpha}+ a_{{\alpha'}}+ \frac{2}{3}( D_{z_{\underline{\alpha}}, j,x} 
+ D_{z_{\underline{\alpha}}, j,y}  + D_{z_{\underline{\alpha}}, j,u} + D_{x,y,u} ) 
-\frac{1}{3} ( D_{z_{\underline{\alpha}}, x,y} 
+ D_{z_{\underline{\alpha}}, y,u}  + D_{z_{\underline{\alpha}}, x,u} + D_{j,x,y} 
+ D_{j,y,u}  + D_{j,x,u} ) =$

\bigskip

$\stackrel{\tiny de\!f.  \; o\!f\; D_{z_{\underline{\alpha}}, \cdot, \cdot} }{=}\;  
 a_{\alpha}+ a_{{\alpha'}} + \frac{2}{3}( D_{\alpha, j,x} 
+ D_{\alpha, j,y}  + D_{\alpha, j,u} + D_{x,y,u} ) -2 a_{\alpha}
-\frac{1}{3} ( D_{\alpha, x,y} 
+ D_{\alpha, y,u}  + D_{\alpha, x,u} + D_{j,x,y} 
+ D_{j,y,u}  + D_{j,x,u}) + a_{\alpha} =$

\bigskip

$\stackrel{\tiny  \;\; Rem. \ref{dt}  }{=}\;
  D_{\alpha, j}+  a_{{\alpha'}}=$

\bigskip

$\stackrel{\tiny  \forall  r,s \in  [n], 
\;\; Rem. \ref{dt}  }{=}\;
 \frac{2}{3}( D_{\alpha, j, \alpha'} 
+ D_{\alpha, j,y}  + D_{\alpha, j,u} + D_{ \alpha',y,u} ) 
-\frac{1}{3} ( D_{\alpha,  \alpha',y} 
+ D_{\alpha, y,u}  + D_{\alpha,  \alpha',u} + D_{j, \alpha',y} 
+ D_{j,y,u}  + D_{j, \alpha',u}) +
\frac{1}{2} ( D_{ \alpha , 
\alpha'}+ D_{ \alpha', r,s} -D_{ \alpha,r,s })=
$

\bigskip

$\stackrel{\tiny  \forall r,s \in  [n],\;\;  Rem. \ref{dt}}{=}\;
 \frac{2}{3}( D_{\alpha, j, \alpha'} 
+ D_{\alpha, j,y}  + D_{\alpha, j,u} + D_{ \alpha',y,u} ) 
-\frac{1}{3} ( D_{\alpha,  \alpha',y} 
+ D_{\alpha, y,u}  + D_{\alpha,  \alpha',u} + D_{j, \alpha',y} 
+ D_{j,y,u}  + D_{j, \alpha',u} )+
 \frac{1}{3}( D_{\alpha,  \alpha', j} 
+ D_{\alpha,  \alpha' ,y}  + D_{\alpha,  \alpha' ,u} + 
D_{ j,y,u} ) 
-\frac{1}{6} ( D_{\alpha,  j,u} 
+ D_{\alpha,j, y}  + D_{\alpha, y,u} + D_{ \alpha',j,u} 
+ D_{ \alpha', j,y}  + D_{\alpha',y,u} )
+
\frac{1}{2} (  D_{ \alpha', r,s} -D_{ \alpha,r,s })
= D_{\alpha , \alpha', j} $

\bigskip

B2) $D_{ \alpha , \alpha' , \beta} (R) 
 \stackrel{\tiny de\!f. \; o\!f\; R}{=}\;  a_{\alpha}+ a_{\alpha'}
 + D_{z_{\underline{\alpha}} , z_{\underline{\beta}}}(T) + a_{\beta}
=$

\bigskip

$\stackrel{\tiny \forall\;x,y,u \in L, \;\; Rem. \ref{dt} }{=}\; 
 a_{\alpha}+ a_{{\alpha'}} + a_{\beta}+ \frac{2}{3}[ D_{z_{\underline{\alpha}}, z_{\underline{\beta}},x}(T) 
+ D_{z_{\underline{\alpha}}, z_{\underline{\beta}},y}(T)  + D_{z_{\underline{\alpha}}, z_{\underline{\beta}},u}(T) + D_{x,y,u} (T)]
-\frac{1}{3} [ D_{z_{\underline{\alpha}}, x,y}(T) 
+ D_{z_{\underline{\alpha}}, y,u}(T)  + D_{z_{\underline{\alpha}}, x,u}(T) + D_{z_{\underline{\beta}},x,y} (T)
+ D_{z_{\underline{\beta}},y,u} (T) + D_{z_{\underline{\beta}},x,u} (T)] $

\bigskip

$\stackrel{\tiny ind. \;ass. }{=}\; 
 a_{\alpha}+ a_{{\alpha'}}+ a_{\beta}+ \frac{2}{3}( D_{z_{\underline{\alpha}}, z_{\underline{\beta}},x} 
+ D_{z_{\underline{\alpha}}, z_{\underline{\beta}},y}  + D_{z_{\underline{\alpha}}, z_{\underline{\beta}},u} + D_{x,y,u} ) 
-\frac{1}{3} ( D_{z_{\underline{\alpha}}, x,y} 
+ D_{z_{\underline{\alpha}}, y,u}  + D_{z_{\underline{\alpha}}, x,u} + D_{z_{\underline{\beta}},x,y} 
+ D_{z_{\underline{\beta}},y,u}  + D_{z_{\underline{\beta}},x,u} )= $

\bigskip

$
\stackrel{\tiny de\!f.\; o\!f\; 
 \; D_{z_{\underline{\alpha}}, \cdot, 
\cdot}  \;and\;  o\!f\;  D_{z_{\underline{\alpha}},z_{\underline{\beta}}  , 
\cdot} }{=}\;  
 a_{\alpha}+ a_{{\alpha'}} + a_{\beta}+ \frac{2}{3}( D_{\alpha, \beta,x} 
+ D_{\alpha, \beta,y}  + D_{\alpha, \beta,u} + D_{x,y,u} ) -2 a_{\alpha}
-2 a_{\beta}
-\frac{1}{3} ( D_{\alpha, x,y} 
+ D_{\alpha, y,u}  + D_{\alpha, x,u} + D_{\beta,x,y} 
+ D_{\beta,y,u}  + D_{\beta,x,u}) + a_{\alpha}  + a_{\beta} =$

\bigskip

$\stackrel{\tiny  Rem. \ref{dt}  }{=}\;
 D_{\alpha, \beta}
+ a_{{\alpha'}}=$

\bigskip

$\stackrel{\tiny  \forall  r,s \in  [n], \;\; Rem. \ref{dt}  }{=}\;
 \frac{2}{3}( D_{\alpha, \beta, \alpha'} 
+ D_{\alpha, \beta,y}  + D_{\alpha, \beta,u} + D_{ \alpha',y,u} ) 
-\frac{1}{3} ( D_{\alpha,  \alpha',y} 
+ D_{\alpha, y,u}  + D_{\alpha,  \alpha',u} + D_{\beta, \alpha',y} 
+ D_{\beta,y,u}  + D_{\beta, \alpha',u}) +
\frac{1}{2} ( D_{ \alpha , 
\alpha'}+ D_{ \alpha' , r,s} -D_{ \alpha,r,s })=
$

\bigskip

$\stackrel{\tiny  \forall r,s \in  [n],\;\; Rem. \ref{dt} }{=}\;
 \frac{2}{3}( D_{\alpha, \beta, \alpha'} 
+ D_{\alpha, \beta,y}  + D_{\alpha, \beta,u} + D_{ \alpha',y,u} ) 
-\frac{1}{3} ( D_{\alpha,  \alpha',y} 
+ D_{\alpha, y,u}  + D_{\alpha,  \alpha',u} + D_{\beta, \alpha',y} 
+ D_{\beta,y,u}  + D_{\beta, \alpha',u} )+
 \frac{1}{3}( D_{\alpha,  \alpha', \beta} 
+ D_{\alpha,  \alpha' ,y}  + D_{\alpha,  \alpha' ,u} + 
D_{ \beta,y,u} ) 
-\frac{1}{6} ( D_{\alpha,  \beta,u} 
+ D_{\alpha,\beta, y}  + D_{\alpha, y,u} + D_{ \alpha',\beta,u} 
+ D_{ \alpha', \beta,y}  + D_{\alpha',y,u} )
+
\frac{1}{2} (  D_{ \alpha' , r,s} -D_{ \alpha,r,s })
=
 D_{\alpha , \alpha', \beta} $

\bigskip

B3) $D_{ \alpha , \alpha' , \alpha''} (R)  
\stackrel{\tiny de\!f. \; o\!f\; R}{=}\; 
a_{\alpha} + a_{\alpha'} + a_{\alpha''} = $

\bigskip

$=\frac{1}{2} 
( D_{\alpha, \alpha'}  + D_{\alpha , x,y} -  D_{\alpha' , x,y} ) +
 \frac{1}{2} 
( D_{\alpha', \alpha''}  + D_{\alpha', x,y} -  D_{\alpha'' , x,y} ) +
 \frac{1}{2} 
( D_{\alpha'', \alpha}  + D_{\alpha'', x,y} -  D_{\alpha, x,y} ) =$

\bigskip
$= \frac{1}{2} ( 
D_{\alpha, \alpha'}  + D_{\alpha' , \alpha''} +  D_{\alpha'' , \alpha} ) 
 =D_{\alpha , \alpha', \alpha''} $, 

where the last equality follows from Remark \ref{dt} with 
$\{r,s,u\} =  \{\alpha'', x,y  \}$ for $D_{\alpha, \alpha'}$,
$\{r,s,u\} =  \{\alpha', x,y  \}$ for $D_{\alpha, \alpha''}$,
$\{r,s,u\} =  \{\alpha, x,y  \}$ for $D_{\alpha', \alpha''}$,
for some $x,y \in [n]$.







\hfill \framebox(7,7)

\begin{remark}
Proposition \ref{starbell} assures us that if we have a positive-weighted
tree $E$ with no  vertices of degree $2$ and 
we consider
the set $\{D_{i,j}(E)\}$  (or $\{ D_{i,j,k}(E) \}$) 
then the tree $R$ coming from  the above theorems is equal to $E$
(or isomorphic to $E$ if we allow vertices of degree $2$ in $E$).
\end{remark}

\section{A modification of Neighbours-Joining method}

We describe shortly  Nei-Saitou's 
Neighbours-Joining algorithm (N-J algorithm for short) 
to reconstruct trees from the data $ \{D_{i,j}\}$ (see 
\cite{NS}, \cite{SK} and \cite{PSt2}).

Let $T$ be a tree with leaves $1,...n$ 
and let $D_{i,j}= D_{i,j} (T)$ for all $i,j \in [n]$.
Define
 $$S_{i,j}:= (n-2) D_{i,j} - \sum_{k \neq i} D_{i,k} - \sum_{k \neq j}
 D_{j,k} $$
for $i \neq j$. 

\begin{lemma} \label{lem}
If $i$ and $j$ are neighbours,  then $S_{i,k}- S_{i,j} \geq 0 $
 for all $k$ and $S_{i,k}- S_{i,j} = 0 $  if and only if
$ k $ is  a neighbour of $i$ and $j$.
\end{lemma}

\begin{theorem}
If  $i$ and $j$ are such that $S_{i,j} $ is a minimum in  the matrix $S$, 
then   $i$ and $j$ are neighbours.
\end{theorem}

The algorithm works in the following way: from $\{D_{i,j}\}$ calculate 
the matrix $S$, defined except on the diagonal;  
let $i,j$ be a couple such that  $S_{i,j} $ is a minimum in  
the matrix $S$; 
by the theorem, we know that then $i$ and $ j$
are neighbours; one can calculate the lengths of their twigs by the formula
$ \frac{1}{2} (D_{i,j}+ D_{i,x} -D_{j,x})$ for any leaf
$x$ and then  consider, instead of 
$ [n]$, the set $ L= ([n] -\{ i,j\})  
\cup \{ z\} $. Define  $D_{z,y} = D_{i,y}-
 \frac{1}{2} (D_{i,j}+ D_{i,y} -D_{j,y})=
\frac{1}{2}(
   -D_{i,j}+ D_{i,y} +D_{j,y})$ for any $y$ 
and iterate the process.

Observe that 
to calculate $S$ we need $O(n^2) $ elementary operations and
to find the minimum of the coefficients of $S$, we need also other 
$O(n^2) $ elementary 
operations. So, to find ONE cherry we need $ O(n^2)$  elementary operations.
(Then we have to repeat  the operations on $L$, whose 
cardinality is $n-1$. So in total we need $O(n^2 + (n-1)^2+.....)= O(n^3) $ 
elementary operations.)

\bigskip

We can modify  the N-J method in the following way.

Given $\{D_{i,j}\} $, calculate the matrix  $S$ 
(defined except on the diagonal).

For every column $S^{(j)}$ of $S$ calculate the minimum 
$m_j$;  this minimum will be attained perhaps more than one 
time.
Choose  the first  index $i$ such that $S_{i,j}$ attains $m_j$, 
precisely let $i_j = min\{i | \; S_{i,j} = m_j\}$.
For these $n$ couples 
of indices $ (i_j,j) $ calculate the difference between the 
$i_j$-th column and the $j$-th column of the matrix $ D $.

- If this difference has all coefficients equal (except from the 
$j$-th and the 
$i_j$-th entries), then  $i_j$ and $j$ are
neighbours by Proposition \ref{starbell}.

- If this difference has not all coefficients equal (except from the 
$j$-th and the $i_j$-th entries), then  $j$ is not in any cherry,
in fact: by Proposition \ref{starbell},  $i_j$ and $j$ are not 
neighbours; besides   if $j$ were in a cherry, let $i$ be one of its 
neighbours, then we would have $ S_{i,j} \leq S_{k,j}$ for any $ k$  by Lemma
\ref{lem}, in particular $S_{i,j} \leq S_{i_j, j} $; thus, since  
 $S_{i_j, j}$ attains the minimum in its column, we must have 
$S_{i,j} = S_{i_j, j} $, but, by Lemma \ref{lem},
 this would imply that $i_j$ would be in the cherry $i,j$, but this is absurd 
because $i_j$ and $j$ are not neighbours. 

So, of the $n$ couples 
of indices $ (i_j,j) $, consider only the ones such that 
 the difference between the 
$i_j$-th column and the $j$-th column of the matrix $ D $
 has all coefficients equal (except from the $j$-th and the 
$i_j$-th entries).

In this way one can easily see that, for every cherry, we have 
found out all its 
leaves, in fact: 
let $i_1,...., i_k \in [n]$ be the leaves of a cherry; suppose 
for instance that $ i_1 <.... < i_k$;   
for the column $i_1$, we will get the 
couple $(i_2,i_1) $ and the difference between the 
$i_2$-th column and the $i_1$-th column of the matrix $ D $
 will have all coefficients equal (except from the $i_1$-th and the 
$i_2$-th entries) so the couple $(i_2,i_1) $ will be taken;
for the column $i_2$ we will get 
the couple $(i_1,i_2) $ and the difference between the 
$i_1$-th column and the $i_2$-th column of the matrix $ D $
 will have all coefficients equal (except from the $i_1$-th and the 
$i_2$-th entries) so the couple $(i_2,i_1) $ will be taken;
for the column $i_3$ the couple $(i_1,i_3) $ 
and the difference between the 
$i_1$-th column and the $i_3$-th column of the matrix $ D $
 will have all coefficients equal (except from the $i_1$-th and the 
$i_3$-th entries) so the couple $(i_1,i_3) $ will be taken;
and so on...

So, by $O(n^2)$ elementary 
operations in all, we have found out ALL the cherries of the tree. 

Once found all the cherries, one can calculate the lengths of all the twigs as
 before, for instance if $i$ and $j$ are neighbours 
the length of the twig of $i$ will be
 $  \frac{1}{2} ( D_{ i , j}+ D_{ i , x} -D_{j ,x})$ for any leaf
$x$ and then substitute, 
for every cherry, to all the element of the cherry
a new point and iterate the process.

\bigskip

Obviously we can allow some ``noise'', by searching, instead of the couples 
$(i,j) $ such that the difference of $i$-th column and the $j$-th column 
of the matrix $ D $  has all coefficients equal, the couples 
$(i,j) $ such that the coefficients of the difference of 
$i$-th column and the $j$-th column of the matrix $ D $ 
differ less than some $\epsilon$. 

\bigskip

Obviously this ``pruning'' version of N-J algorithm can be a bit quicker 
than the original one  only in the cases where at every step there are a lot 
of cherries; besides it could be useful when we can't calculate all the tree but 
we want to know at least the list of all the cherries. Finally we observe that 
perhaps an underestimate of $D_{i,j}$ for some $i,j$ can be misleading 
in the reconstuction of the tree by classical N-J algorithm and that the 
characterization of neighbours given by Proposition \ref{starbell} and the 
above described version of N-J algorithm can 
detect such misleading.

\bigskip

In \cite{LYP}, \cite{PSt} and \cite{PSt2} a generalization of N-J algorithm 
to get trees from $k$-weights is described; namely let 
$$ S_{i,j}= \frac{n-2}{k-1} \sum_{\tiny
\begin{array}{c}
R \; (k-2)\!-\!subset \\
 \; o\!f \; [n] -\{i,j\}
\end{array}}
 D(i,j,R) - \sum_{\tiny
\begin{array}{c}
R \; (k-1)\!-\!subset 
\\ \; o\!f \; [n] -\{i\}
\end{array}
}
 D(i,R) - \sum_{\tiny
\begin{array}{c}
R \; (k-1)\!-\!subset \\ 
\; o\!f \; [n] -\{j\}
 \end{array}}
D(j,R)$$

\smallskip

If $i,j$ are such that $S_{i,j}$ is the minimum in the matrix $S$ 
then $i$ and $j$ are neighbours (to calculate $S$ and then to find out 
a cherry, we need $O(n^k)$ elementary operations) and we can proceed as before.

Again by using Proposition \ref{starbell} we can modify the algorithm 
 finding  all the cherries by $O(n^k)$ elementary operations.

\bigskip

\end{document}